%Paper: funct-an/9405001
%From: exel@ime.usp.br
%Date: Mon, 2 May 94 18:05:52-030

  % This is a Plain TeX file
  % FONTS
  \newcount\fontset
  \fontset=2
  \def\dualfont#1#2#3{\font#1=\ifnum\fontset=1 #2\else#3\fi}
  \dualfont\eightrm {cmr8} {cmr7}
  \dualfont\eightit {cmti8} {cmti10}
  \dualfont\tensc {cmcsc10} {cmcsc10}
  \dualfont\titlefont {cmss12} {cmr17}

  % GLOBAL SETTINGS
  \magnification=\magstep1
  \nopagenumbers
  \voffset=2\baselineskip
  \advance\vsize by -\voffset
  \headline{\ifnum\pageno=1 \hfil
    \else\ifodd\pageno\tensc\hfil twisted partial actions \hfil\folio
    \else\tensc\folio\hfil ruy exel\hfil\fi\fi}

  % STANDARD DEFINITIONS
  \def\<{\langle}
  \def\>{\rangle}
  \def\:{\colon}
  \def\*{\otimes}
  \def\+{\oplus}
  \def\x{\times}
  \def\c{\subseteq}
  \def\arw{\rightarrow}
  \def\cstar{$C^*$}
  \def\square{\hbox{$\sqcap\!\!\!\!\sqcup$}}

  % CONTROL SEQUENCES
  \def\items{\begingroup\parskip 5pt\parindent 20pt \mathsurround=3pt}
  \def\itm#1{\item{#1)}}
  \def\vg#1{\ifx#1*\empty\else\ifx#1,,\else \ifx#1..\else
    \ifx#1;;\else \ifx#1::\else \ifx#1''\else \ifx#1--\else
    \ifx#1))\else \ifx#1/\else { }#1\fi\fi\fi \fi\fi\fi\fi\fi\fi}
  \newcount\bibno \bibno=0
  \def\newbib#1{\advance\bibno by 1 \edef#1{\number\bibno}}
  \def\stdbib#1#2#3#4#5#6#7#8{\smallskip \item{[#1]} #2, ``#3'',
    {\sl#4} {\bf#5} (#6), #7--#8.}
  \def\bib#1#2#3#4{\smallskip \item{[#1]} #2, ``#3'', {#4.}}
  \def\cite#1{{\rm[\bf #1\rm]}}
  \def\scite#1#2{{\rm[\bf #1\rm, #2]}}
  \def\lcite#1{(#1)}
  \newcount\secno \secno=0
  \newcount\stno \stno=0
  \outer\def\beginsection#1\par{
    \stno=0
    \advance\secno by 1
    \vskip0pt plus.1\vsize\penalty-250
    \vskip0pt plus-.1\vsize\bigskip\vskip\parskip
    \message{\number\secno.#1}\centerline{\tensc\number\secno. #1}
    \nobreak\bigskip}
  \def\state#1 #2\par{\advance\stno by 1\medbreak\noindent
    {\bf\number\secno.\number\stno.\enspace #1.\enspace}{\sl
    #2}\medbreak}
  \def\nstate#1 #2 #3\par{\state{#2} {#3}\par
    \edef#1{\number\secno.\number\stno}}
  \def\beginproof{\medbreak\noindent{\it Proof.\enspace}}
  \def\endproof{\ifmmode\eqno\square\else\hfill\square\medbreak\fi}

  % MISCELLANEOUS DEFINITIONS

  \def\piso{partial isometry\vg}
  \def\spb{semidirect product bundle\vg}
  \def\pisos{partial isometries\vg}
  \def\lc{locally compact\vg}
  \def\tpa{twisted partial action\vg}
  \def\tpas{\tpa*s\vg}
  \def\ab{algebraic bundle\vg}
  \def\abs{\ab*s\vg}
  \def\tro{TRO\vg}
  \def\tros{\tro*s\vg}

  \def\B{{\cal B}}
  \def\D{{\cal D}}
  \def\L{{\cal L}}
  \def\R{{\cal R}}
  
  \def\K{{\cal K}}
  \def\M{{\cal M}}
  \def\E{{\cal E}}

  \def\d{\delta}
  \def\e{\varepsilon}
  \def\g{\gamma}
  \def\t{\theta}
  \def\Rho{\hbox{\rm P}}
  \def\rhodag{\rho^\dagger}
  \def\lamdag{\lambda^\dagger}

  \def\assoc{\sim}
  \def\sassoc{\buildrel s \over \assoc}
  \def\ltwo{l_2}
  \def\Link{{\rm Link}}
  \def\cob{{\cal C}_0(\B)}
  \def\cod{{\cal C}_0(\D)}
  \def\codinv{{\cal C}_0(\D^{-1})}

  \def\({\left(}
  \def\){\right)}
  \def\for#1{,\quad #1}
  \def\and{\quad\hbox{and}\quad}
  \def\rr{(resp.~right)\vg}
  \def\mat#1#2#3#4{\pmatrix{#1 & #2 \cr #3 & #4}}
  \def\tri#1{|\kern-1 pt|\kern-1 pt|#1|\kern-1 pt|\kern-1 pt|}
  \def\m#1#2#3{#1\cdot#2\cdot#3}
  \def\stackrel#1#2{\buildrel #1 \over #2}

  % REFERENCES
  %
  \newbib{\BGR}
  \newbib{\BMS}
  \newbib{\Newpim}
  \newbib{\Fell}
  \newbib{\McCl}
  \newbib{\PR}
  \newbib{\Pedersen}
  \newbib{\Quigg}
  \newbib{\Induc}
  \newbib{\MoritaOne}
  \newbib{\MoritaTwo}
  \newbib{\ZM}
  \newbib{\Zettl}

  % START OF TEXT
  \null
  \begingroup
  \def\c{\centerline}

  \titlefont \baselineskip=3ex

  % TITLE
  %
  \c{Twisted Partial Actions}
  \c{A Classification of Stable C*-Algebraic Bundles}

  \bigskip
  \c{\tt Peliminary Version}

  \bigskip \tensc \baselineskip=3ex

  \c{Ruy Exel}

  \bigskip \eightit \baselineskip=3ex

  \c{Departamento de Matem\'atica}
  \c{Universidade de S\~ao Paulo}
  \c{C.~P.~ 20570}
  \c{01452-990 S\~ao Paulo Brazil}
  \c{e-mail: exel@ime.usp.br}

  % ABSTRACT
  %
  \bigskip\bigskip\eightrm\baselineskip=3ex
  \midinsert \narrower
  We introduce the notion of continuous \tpas of a \lc group on a
C*-algebra. With such, we construct an associated C*-\ab called the
\spb.  Our main theorem shows that, given any C*-\ab which is second
countable and whose unit fiber algebra is stable, there is a
continuous \tpa of the base group on the unit fiber algebra, whose
associated \spb is isomorphic to the given one.
  \endinsert \endgroup

  \beginsection Introduction

  A \cstar-\ab is, roughly speaking, a natural generalization of the
concept of graded \cstar-algebras, to the case when the grading group
is a \lc group. A \cstar-\ab\null $\B$ over the group $G$ consists,
therefore, of a collection of Banach spaces $\(B_t\)_{t\in G}$ which
are glued together to form a Banach bundle \scite{\Fell}{II.13.4} and
which moreover comes equipped with a family of multiplication
operations
  $$\cdot\ \:B_r\x B_s \arw B_{rs} \for r,s \in G$$
  and a family of involution operations
  $$*\ \:B_t \arw B_{t^{-1}} \for t \in G$$
  all of which are continuous with respect to $r$, $s$ and $t$,
satisfying axioms that are modeled after the properties which would be
satisfied, were the $B_t$'s the grading subspaces of a graded
*-algebra $\bigoplus_{t\in G}B_t$.

  \cstar-\abs naturally occur in a large number of situations. First
of all, they show up in connection with the theory of group
representations, as carefully described in the comprehensive two
volume book by Fell and Doran, listed below as reference number
\cite{\Fell}.
  In the discrete group case, there is a very close relationship
(although not a perfect equivalence \scite{\Fell}{VIII.16.12}) between
\cstar-\abs and graded \cstar-algebras. The latter, in turn, appears
in connection to the theory of group actions on \cstar-algebras. In
fact, whenever a compact abelian group $K$ acts on a \cstar-algebra
$B$, then there is a natural grading on $B$ by the dual group
$G=\hat{K}$, given by the spectral subspaces (see, for example
\scite{\Newpim}{Section 2} for the case of the circle group). More
generally, Quigg \cite{\Quigg} has shown that a co-action of a
discrete group (a concept which generalizes actions of the not always
visible compact dual group) also yields graded \cstar-algebras.

  In this work we propose to extend our earlier work \cite{\Newpim} on
circle actions (i.e. {\bf Z}-graded algebras) to the context of
general \cstar-\abs. As a result, we obtain a classification theorem
which exhibits any \cstar-\ab, satisfying certain mild hypothesis, as
the {\it \spb\null} for a continuous {\it\tpa\null} of the base group
on the unit fiber algebra.  The concept of \tpas is a simultaneous
generalization of the twisted crossed-products of Zeller-Meier
\cite{\ZM} on one hand and the partial actions which we introduced in
\cite{\Newpim} on the other
   (see also the work of Packer and Raeburn \cite{\PR} for the twisted
case, and McClanahan's work \cite{\McCl} on partial actions of
discrete groups).

  A continuous \tpa of a \lc group $G$ on a \cstar-algebra $A$
consists of a family $\{D_t\}_{t\in G}$ of closed two sided ideals of
$A$, a family $\{\t_t\}_{t\in G}$ of isomorphisms from $D_{t^{-1}}$ to
$D_t$ and a ``cocycle'' $w=\{w(r,s)\}_{(r,s)\in G\x G}$, where each
$w(r,s)$ is a unitary multiplier of the ideal $D_r\cap D_{rs}$, which
satisfy properties similar to the axioms of twisted actions (see below
for more details).

  Given such an object, we construct a \cstar-\ab called the \spb of
$A$ and $G$, generalizing \scite{\Fell}{VIII.4}. The power of this
construction is such that we are able, in turn, to show, in our main
Theorem, that every second countable \cstar-\ab, whose unit fiber
algebra is stable, can be obtained as the result of our construction.

  The requirement that the unit fiber algebra be compact can obviously
be dropped if one is willing to ``stabilize'' the given bundle, by
tensoring it with the algebra of compact operators on a separable
Hilbert space. Given this, our theorem can then be applied to
virtually all \cstar-\abs.

  As far as introducing a generalized partial crossed-product algebra
(which we did in \cite{\Newpim} in the case of partial actions of the
integers or was done in \cite{\McCl} for partial actions by general
discrete groups), observe that the usual process of forming the
crossed-product of a \cstar-algebra by a group action \cite{\Pedersen}
can be divided in two steps, the first one being the construction of
the associated \spb as in \scite{\Fell}{VIII.4}. The second step is
then to form the cross-sectional algebra \scite{\Fell}{VII.5}, which
is a process that can be applied to any \cstar-\ab, irrespective of
how it came about.

  That is, we may stop short of defining the concept of
crossed-products by \tpas, since the associated \spb can then be fed
to the machinery of cross-sectional algebras, which would, via a
standard procedure, produce what we would call the crossed-product
\cstar-algebra of a \lc group by a continuous \tpa. In particular, we
completely avoid, in this way, the usual problems caused by
non-amenable groups.

  \beginsection The Discrete Group case

  The axioms for continuous \tpas of \lc groups and the basic work
leading to the construction of the associated \cstar-\ab can be
divided in two distinct parts, the first one relating to algebraic
properties and the second referring to topological aspects.
  In order to organize the exposition, we have, therefore, chosen to
break up the presentation of the definition and basic properties of
this concept in two sections, the present one being dedicated to the
algebraic considerations. We thus restrict our initial discussion to
groups without topology or, what amounts to the same, to discrete
groups.

  Let $A$ be a \cstar-algebra and let $G$ be a discrete group.

  \nstate{\DefTPA} Definition A \tpa of $G$ on $A$ is a triple
  $$\Theta = \(\{D_t\}_{t\in G}, \{\t_t\}_{t\in G},
\{w(r,s)\}_{(r,s)\in G\x G}\)$$
  where for each $t$ in $G$, $D_t$ is a closed two sided ideal in $A$,
$\t_t$ is a *-isomorphism from $D_{t^{-1}}$ onto $D_t$ and for each
$(r,s)$ in $G\x G$, $w(r,s)$ is a unitary multiplier of $D_r\cap
D_{rs}$, satisfying the following postulates, for all $r$, $s$ and $t$
in $G$
  \items
  \itm{a} $D_e = A$ and $\t_e$ is the identity automorphism of $A$.
  \itm{b} $\t_r(D_{r^{-1}}\cap D_s) = D_r\cap D_{rs}$
  \itm{c} $\t_r(\t_s(a)) = w(r,s)\t_{rs}(a)w(r,s)^* \for{a\in
D_{s^{-1}} \cap D_{s^{-1}r^{-1}}}$
  \itm{d} $w(e,t)=w(t,e)=1$
  \itm{e} $\t_r(aw(s,t))w(r,st)=\t_r(a)w(r,s)w(rs,t) \for{a \in
D_{r^{-1}}\cap D_s\cap D_{st}}$
  \endgroup

  Our goal is to construct, given a \tpa of $G$ on $A$, a \cstar-\ab
over $G$. We recall below, the definition of this concept in the
special case of discrete groups \scite{\Fell}{II.13.1, II.13.4,
VIII.2.2, VIII.3.1, VIII.16.2}.

  \nstate{\CstarAlgBundle} Definition A \cstar-\ab over a discrete
group $G$ is a collection of Banach spaces $\{B_t\}_{t\in G}$ together
with a multiplication operation
  $$\cdot\ \:\B\x \B \arw \B$$ and an involution
  $$*\ \:\B\arw \B$$
  where $\B$ is the disjoint union of all $B_t$'s, satisfying for all
$r$, $s$ and $t$ in $G$ and $b$ and $c$ in $\B$
  \items
  \itm{i}$B_rB_s\c B_{rs}$
  \itm{ii} The product $\cdot$ is bilinear on $B_r\x B_s$ to $B_{rs}$.
  \itm{iii} The product on $\B$ is associative.
  \itm{iv} $\|bc\| \leq \|b\| \|c\|$
  \itm{v} $(B_t)^* \c B_{t^{-1}}$
  \itm{vi} * is conjugate-linear from $B_t$ to $B_{t^{-1}}$
  \itm{vii} $(bc)^* = c^* b^*$
  \itm{viii} $b^{**}=b$
  \itm{ix} $\|b^*\| = \|b\|$
  \itm{x} $\|bb^*\| = \|b\|^2$
  \itm{xi} $bb^* \geq 0$ in $B_e$
  \endgroup

  Axioms (i)--(iv) define what is called a Banach \ab\null
\scite{\Fell}{VIII.2.2}. Adding (v)--(ix) gives the definition of a
Banach *-\ab\ \scite{\Fell}{VIII.3.1} while the last two properties
characterize \cstar-\abs\null\scite{\Fell}{VIII.16.2}. In the above
definition we have omitted all references to continuity, since we are,
for the time being, considering exclusively discrete groups. See
section \lcite{3} below for the general case.

  We shall denote both the family $\{B_t\}_{t\in G}$ and the disjoint
union of the $B_t's$ by $\B$ as this will not bring any confusion.

  Fix, from now on, a \cstar-algebra $A$, a discrete group $G$ and a
\tpa of $G$ on $A$ given by
  $$\Theta = \(\{D_t\}_{t\in G}, \{\t_t\}_{t\in G},
\{w(r,s)\}_{(r,s)\in G\x G}\).$$
  As a first step in constructing a \cstar-\ab from $\Theta$, let
  $$\B=\{(a,s)\in A\x G: a\in D_s\},$$
  and for each $t$ in $G$, let $B_t$ be the subset of $\B$ formed by
all the $(a,s)$ with $s=t$. We will also use the notation $a\d_t$ for
$(a,t)$, whence $B_t = D_t\d_t$.

  There is an obvious bijection between $B_t$ and $D_t$, through which
we give $B_t$ the structure of a Banach space.

  Let us define the multiplication operation on $\B$ by
  $$(a_r\d_r)*(b_s\d_s) = \t_r\(\t_r^{-1}(a_r)b_s\)w(r,s)\d_{rs}$$
  for $a_r$ in $D_r$ and $b_s$ in $D_s$. It is important to remark
that the term
  $\t_r\(\t_r^{-1}(a_r)b_s\)$ belongs to $D_r\cap D_{rs}$ by
\lcite{\DefTPA.b} and hence that multiplication of this term by
$w(r,s)$ is well defined. The result lying again in $D_r\cap D_{rs}$,
guarantees that the right hand side above in fact gives an element in
$B_{rs}$. We have thus verified \lcite{\CstarAlgBundle.i}.

  The involution on $\B$ is defined by
  $$(a_t\d_t)^* = \t_t^{-1}(a_t^*)w(t^{-1},t)^*\d_{t^{-1}}$$
  for $a_t$ in $D_t$. It should be noted that this gives $(B_t)^* \c
B_{t^{-1}}$, hence proving \lcite{\CstarAlgBundle.v}.

  Also observe that \lcite{\DefTPA.c} with $r=t^{-1}$, $s=t$ and
$a=\t_t^{-1}(a_t^*)$ provides
  $$\t_{t^{-1}}(a_t^*)=w(t^{-1},t) \t_t^{-1} (a_t^*) w(t^{-1}, t)^*$$
  and so the definition above is equivalent to
  $$(a_t\d_t)^*=w(t^{-1},t)^* \t_{t^{-1}} (a_t^*)\d_{t^{-1}}.$$

  One should be careful not to mistake $\t_{t^{-1}}$ for $\t_t^{-1}$
which, as seen above, coincide only up to conjugation by
$w(t^{-1},t)$.

  This specifies all of the required ingredients of a \cstar-\ab and
we must therefore verify the validity of properties (i)--(xi)
above. Apart from (i) and (v) which have already been checked, note
that (ii), (iv), (vi) and (ix) can all be proved without much effort.
In contrast, proving associativity and the anti-multiplicativity of
the involution is a bit of a challenge.

  \state Lemma If \/ $\{u_i\}_i$ is an approximate identity (always
assumed to be self-adjoint and of norm one) for $D_{r^{-1}}$ and if
$a_r$ and $b_s$ are elements of $D_r$ and $D_s$, respectively, then
  $$(a_r\d_r)*(b_s\d_s) = \lim_ia_r\t_r(u_ib_s)w(r,s)\d_{rs}.$$

  \beginproof We have
  $$(a_r\d_r)*(b_s\d_s) = \t_r\(\t_r^{-1}(a_r)b_s\)w(r,s)\d_{rs}$$
  $$ =\lim_i\t_r\(\t_r^{-1} (a_r)u_ib_s\)w(r,s)\d_{rs}
  = \lim_ia_r\t_r(u_ib_s)w(r,s)\d_{rs}.\endproof$$

  \state Proposition The multiplication defined above is associative.

  \beginproof Let $a_r$, $b_s$ and $c_t$ be in $D_r$, $D_s$ and $D_t$,
respectively. Also let $\{u_i\}_i$ be an approximate identity for
$D_{r^{-1}}$. Then we have
  $$(a_r\d_r * b_s\d_s) * c_t\d_t
  = \(\lim_ia_r \t_r(u_ib_s)w(r,s)\d_{rs}\)*c_t\d_t = \ldots$$
  Let $x_i=a_r \t_r(u_ib_s)w(r,s)$. So the above equals
  $$\ldots = \lim_ix_i\d_{rs} * c_t\d_t
  = \lim_i\t_{rs}\(\t_{rs}^{-1}(x_i)c_t\)w(rs,t)\d_{rst}=\ldots$$
  Note that $x_i$ is in $D_r\cap D_{rs}$ so that $\t_{rs}^{-1}(x_i)$
is in
  $D_{s^{-1}} \cap D_{s^{-1}r^{-1}}$.
  Let $v_j$ be an approximate identity for
  $D_{s^{-1}} \cap D_{s^{-1}r^{-1}}$.
  So the above equals
  $$\ldots=\lim_i \lim_j \t_{rs}\(\t_{rs}^{-1}(x_i)v_jc_t\)w(rs,t)
\d_{rst}$$
  $$= \lim_i \lim_j a_r \t_r(u_ib_s)w(r,s) \t_{rs}(v_jc_t)w(rs,t)
\d_{rst}$$
  $$= \lim_i \lim_j a_r \t_r(u_ib_s)\t_r\(\t_s(v_jc_t)\)w(r,s)
w(rs,t)\d_{rst}$$
  $$ = \lim_i \lim_j a_r \t_r\left[u_ib_s\t_s(v_jc_t)
\right]w(r,s)w(rs,t)\d_{rst}$$
  $$= \lim_i \lim_j a_r \t_r\left[\t_s\(\t_s^{-1}(u_ib_s)v_jc_t\)
\right]w(r,s)w(rs,t) \d_{rst}=\ldots$$
  Note that $u_ib_s$ is in $D_{r^{-1}}\cap D_s$ so that
$\t_s^{-1}(u_ib_s)$ is in $D_{s^{-1}}\cap D_{s^{-1}r^{-1}}$ so the
above equals
  $$\ldots=\lim_i a_r \t_r\left[\t_s\(\t_s^{-1}(u_ib_s)c_t\)
\right]w(r,s)w(rs,t) \d_{rst}=\ldots$$
  Without loss of generality we may assume that $b_s=b_s'b_s''$ where
both $b_s'$ and $b_s''$ belong to $D_s$. So the above equals
  $$\ldots=\lim_i\t_r\left[\t_r^{-1}(a_r)\t_s\(\t_s^{-1}(u_ib_s')\t_s^{-1}
(b_s'')c_t\)\right]w(r,s)w(rs,t)\d_{rst}$$
  $$=\lim_i\t_r\left[\t_r^{-1}(a_r)u_ib_s'\t_s\(\t_s^{-1}(b_s'')c_t\)
\right]w(r,s)w(rs,t)\d_{rst}$$
  $$=\t_r\left[\t_r^{-1}(a_r)b_s'\t_s\(\t_s^{-1}(b_s'')c_t\)
\right]w(r,s)w(rs,t)\d_{rst}$$
  $$=\t_r\left[\t_r^{-1}(a_r)\t_s\(\t_s^{-1}(b_s)c_t\)
\right]w(r,s)w(rs,t)\d_{rst}.$$

  On the other hand
  $$a_r\d_r*(b_s\d_s*c_t\d_t)
  = a_r\d_r*\t_s\(\t_s^{-1}(b_s)c_t\)w(s,t)\d_{st}$$
  $$= \t_r\left[\t_r^{-1}(a_r)\t_s\(\t_s^{-1}(b_s)c_t\)w(s,t)\right]
w(r,st)\d_{rst}=\ldots$$
  If $x=\t_r^{-1}(a_r)\t_s\(\t_s^{-1}(b_s)c_t\)$, then $x$ is in
$D_{r^{-1}}\cap D_s\cap D_{st}$ so by \lcite{\DefTPA.e} the above
equals
  $$\ldots=\t_r(x)w(r,s)w(rs,t)\d_{rst}
  =\t_r\left[\t_r^{-1}(a_r)\t_s\(\t_s^{-1}(b_s)c_t\)
\right]w(r,s)w(rs,t)\d_{rst}.\endproof$$

  Let us now state an identity to be used in the proof of the
anti-multiplicativity of the involution. The proof is omitted as this
is basically a rewriting of \lcite{\DefTPA.e}.

  \nstate{\Aux} Lemma If $f$, $g$ and $h$ are in $G$, then
  $$\t_f\(aw(g,h)^*\) = \t_f(a)w(f,gh)w(fg,h)^*w(f,g)^* \for{a \in
D_{f^{-1}}\cap D_g\cap D_{gh}}.$$

  \state Proposition The involution defined above is
anti-multiplicative.

  \beginproof Let $a_r$ and $b_s$ be in $D_r$ and $D_s$, respectively.
We want to prove that
  $$(a_r\d_r * b_s\d_s)^* = (b_s\d_s)^* * (a_r\d_r)^*.$$
  The left hand side equals
  $$\left[\t_r\(\t_r^{-1}(a_r)b_s\)w(r,s)\d_{rs}\right]^*
  = \t_{rs}^{-1}\left[w(r,s)^*\t_r\(b_s^*\t_r^{-1}(a_r^*)\)\right]
w(s^{-1}r^{-1},rs)^*\d_{s^{-1}r^{-1}}.$$
  Let $x=b_s^*\t_r^{-1}(a_r^*)$. Then $x$ is in $D_s\cap D_{r^{-1}}$
whence $\t_s^{-1}(x)$ is in $D_{s^{-1}}\cap D_{s^{-1}r^{-1}}$ and we
have by axiom \lcite{\DefTPA.c}
  $$\t_r(x)=\t_r\(\t_s\(\t_s^{-1}(x)\)\)
  = w(r,s)\t_{rs}\(\t_s^{-1}(x)\)w(r,s)^*$$
  and thus
  $$(a_r\d_r * b_s\d_s)^* = \t_{rs}^{-1}\left[\t_{rs}\(\t_s^{-1}(x)\)
w(r,s)^*\right]w(s^{-1}r^{-1},rs)^*\d_{s^{-1}r^{-1}}$$
  $$=w(s^{-1}r^{-1},rs)^*\t_{s^{-1}r^{-1}}\left[\t_{rs}\(\t_s^{-1}(x)\)
w(r,s)^*\right]\d_{s^{-1}r^{-1}}=\ldots$$
  Using \lcite{\Aux} with $f=s^{-1}r^{-1}$, $g=r$, $h=s$ and
$a=\t_{rs}\(\t_s^{-1}(x)\)$ the above equals
  $$\ldots=w(s^{-1}r^{-1},rs)^*\t_{s^{-1}r^{-1}}\left[\t_{rs}\(\t_s^{-1}
(x)\)\right] w(s^{-1}r^{-1},rs)w(s^{-1},s)^*w(s^{-1}r^{-1},r)^*
\d_{s^{-1}r^{-1}}$$
  $$=
\t_{rs}^{-1}\left[\t_{rs}\(\t_s^{-1}(x)\)\right]w(s^{-1},s)^*w(s^{-1}
r^{-1},r)^* \d_{s^{-1}r^{-1}}$$
  $$=\t_s^{-1}(x)w(s^{-1},s)^*w(s^{-1} r^{-1},r)^* \d_{s^{-1}r^{-1}}$$
  $$=\t_s^{-1}\(b_s^*\t_r^{-1}(a_r^*)\)w(s^{-1},s)^*w(s^{-1}
r^{-1},r)^* \d_{s^{-1}r^{-1}}.$$

  On the other hand
  $$(b_s\d_s)^* * (a_r\d_r)^* =
\(\t_s^{-1}(b_s^*)w(s^{-1},s)^*\d_{s^{-1}}\) *
\(\t_r^{-1}(a_r^*)w(r^{-1},r)^*\d_{r^{-1}}\)$$
  $$=\(w(s^{-1},s)^*\t_{s^{-1}}(b_s^*)\d_{s^{-1}}\) *
\(\t_r^{-1}(a_r^*)w(r^{-1},r)^*\d_{r^{-1}}\)$$
  $$=\t_{s^{-1}}\left\{\t_{s^{-1}}^{-1}\left[ w(s^{-1},s)^*
\t_{s^{-1}}(b_s^*)\right]\t_r^{-1}(a_r^*)w(r^{-1},r)^*\right\}
w(s^{-1}r^{-1})\d_{s^{-1},r^{-1}}=\ldots$$
  Let $x=\t_{s^{-1}}(b_s^*)$ and $y=\t_r^{-1}(a_r^*)w(r^{-1},r)^*$ and
let $u_i$ be an approximate identity for $D_s$. So the above equals
  $$\ldots=\lim_i
\t_{s^{-1}}\left\{\t_{s^{-1}}^{-1}\left[w(s^{-1},s)^*x\right]u_i
y\right\}w(s^{-1},r^{-1})\d_{s^{-1}r^{-1}}$$
  $$= \lim_i
\t_{s^{-1}}\left\{\t_{s^{-1}}^{-1}\left[w(s^{-1},s)^*x\t_{s^{-1}}
(u_iy) \right]\right\}w(s^{-1},r^{-1})\d_{s^{-1}r^{-1}}$$
  $$=\lim_iw(s^{-1},s)^*x\t_{s^{-1}}(u_iy)w(s^{-1},r^{-1})\d_{s^{-1}r^{-1}}$$
  $$=\lim_iw(s^{-1},s)^*\t_{s^{-1}}(b_s^*u_iy)w(s^{-1},r^{-1})
\d_{s^{-1}r^{-1}}$$
  $$=w(s^{-1},s)^*\t_{s^{-1}}(b_s^*y)w(s^{-1},r^{-1})\d_{s^{-1}r^{-1}}$$
  $$=w(s^{-1},s)^*\t_{s^{-1}}\left[b_s^* \t_r^{-1}(a_r^*)w(r^{-1},r)^*
\right]w(s^{-1},r^{-1})\d_{s^{-1}r^{-1}}=\ldots$$
  Let us now use lemma \lcite{\Aux} once more with $f=s^{-1}$,
$g=r^{-1}$, $h=r$ and $a=b_s^*\t_r^{-1}(a_r^*)$ to conclude that the
above equals
  $$\ldots=w(s^{-1},s)^* \t_{s^{-1}}\(b_s^*\t_r^{-1}(a_r^*)\)
w(s^{-1},e) w(s^{-1}r^{-1},r)^*w(s^{-1},r^{-1})^*
w(s^{-1},r^{-1})\d_{s^{-1}r^{-1}}$$
  $$=\t_s^{-1}\(b_s^*\t_r^{-1}(a_r^*)\)w(s^{-1},s)w(s^{-1}r^{-1},r)^*
\d_{s^{-1}r^{-1}}.\endproof$$

  We leave for the reader to verify the remaining properties, i.e.,
(viii), (x) and (ix) of Definition \lcite{\CstarAlgBundle}. Once this
is done we have proven the main result of this section:

  \state Theorem Given a {\tpa\null}
  $$\Theta = \(\{D_t\}_{t\in G}, \{\t_t\}_{t\in G},
\{w(r,s)\}_{(r,s)\in G\x G}\)$$
  of the discrete group $G$ on $A$, the bundle $\B=\{D_t\d_t\}_{t\in
G}$ is a \cstar-\ab over $G$ with the operations
  $$(a_r\d_r)*(b_s\d_s) = \t_r\(\t_r^{-1}(a_r)b_s\)w(r,s)\d_{rs}$$ and
  $$(a_t\d_t)^* = \t_t^{-1}(a_t^*)w(t^{-1},t)^*\d_{t^{-1}}.$$

  \state Definition The \cstar-\ab constructed above will be called
the \spb of $A$ and $G$ (after \scite{\Fell}{VIII.4}).

  This concludes the algebraic part of our construction.

  \beginsection The Continuous Group case

  From now on we will let $G$ be a \lc topological group.  Of course
one would like to add extra requirements to the definition of \tpas to
account for the topology of $G$. In other words we would like to
require that \tpas be continuous in a suitable sense to be made
precise below.

  To begin with, let us establish the relevant concept of continuity
for a family of subspaces of a given Banach space. Let, therefore, $E$
be a Banach space and let
  $\{E_x\}_{x\in X}$
  be a family of linear subspaces of $E$, indexed by a topological
space $X$.

  \nstate{\DefContSubspace} Definition We say that
  $\{E_x\}_{x\in X}$ is continuous if, for any open set $U \c E$, the
set
  $$\{x\in X: E_x\cap U\neq \emptyset\}$$
  is open in $X$.

  Consider the subset $\E$ of $E\x X$ defined by
  $$\E = \{(v,x)\in E\x X: v\in E_x\}$$
  equipped with the relative topology from $E\x X$. If continuity is
assumed, we claim that $\E$ is a Banach bundle over $X$, as defined in
\scite{\Fell}{II.13.4}. In fact, the properties (i)--(iv) in that
definition are automatically verified as they all follow from the
corresponding facts which hold for the trivial Banach bundle $E\x X$.
The main question hinges on the openness of the bundle projection.

  \nstate{\BundleVsContinuity} Proposition The family $\{E_x\}_{x\in
X}$ is continuous if and only if the bundle projection
  $$\pi\:(v,x)\in \E \arw x\in X$$
  is an open map. In this case $\E$ is a Banach bundle over $X$.

  \beginproof The proof is elementary and hence omitted.\endproof

  Observe that a section $\gamma\:X\arw E$ must necessarily have the
form $\gamma(x) = (\beta(x),x)$ for some function $\beta\:X\arw E$
such that $\beta(x)\in E_x$ for all $x$ in $X$. Since $\E$ has the
relative topology, one sees that $\gamma$ is continuous if and only if
$\beta$ is. Given the very close relationship between $\gamma$ and
$\beta$, they will be deliberately confused with each other.

  \nstate{\LotsSections} Proposition If for any $x_0$ in $X$, any
$\varepsilon>0$ and any $v$ in $E_{x_0}$ there exists a continuous
section $\beta$ such that $\|\beta(x_0)-v\|<\varepsilon$, then
  $\{E_x\}_{x\in X}$
  is continuous. In addition, if $X$ is either \lc or paracompact,
then the converse holds even strongly in the sense that a section
$\beta$ can always be found with $\beta(x_0)=v$.

  \beginproof Let $U$ be an open subset of $E$ and suppose that $x_0$
is such that
  $E_{x_0}\cap U \neq \emptyset$. Pick a continuous section $\beta$
such that $\beta(x_0) \in E_{x_0}\cap U.$ Note that for any $x$ in
$\beta^{-1}(U)$ one has that $\beta(x)\in E_x\cap U$ and consequently
$E_x\cap U$ is non-empty. If we note that $\beta^{-1}(U)$ is open, we
see that the the first part of the statement is proven. Conversely, if
$\{E_x\}_{x\in X}$ is continuous, then by \lcite{\BundleVsContinuity},
$\E$ is a Banach bundle over $X$. The conclusion, then follows from a
result of Douady and Dal Soglio-H\'erault, stating that Banach bundles
over \lc or paracompact base spaces have plenty of continuous sections
\scite{\Fell}{II.13.19}. \endproof

  In connection to this let us define, for future reference, the
concept of pointwise-dense set of sections.

  \state Definition A set $\Gamma$ of continuous sections of a given
Banach bundle is said to be pointwise-dense if for any $x$ in the base
space, the set $\{\gamma(x)\:\gamma\in \Gamma\}$ is dense in the fiber
over $x$.

  Returning to the case of our \tpa, let us assume henceforth that the
collection of ideals $\{D_t\}_{t\in G}$ is continuous in the sense
above. Since the group is assumed to be \lc we will thus have
continuous sections in abundance.

  As in the previous section, let $\B=\{(a,t)\in A\x G: a\in D_t\}$
which we now consider as a topological subspace of $A\x G$. It then
follows from \lcite{\BundleVsContinuity} that $\B$ is a Banach bundle
over $G$.

  Because the inversion map is continuous on $G$, the same reasoning
above shows that $\B^{-1}= \{(a,t)\in A\x G: a\in D_{t^{-1}}\}$ is
also a Banach bundle, and then the family of isomorphisms
  $\{\t_t\}_{t\in G}$ can be used to define a bundle map
$\t\:\B^{-1}\arw \B$.

  \nstate{\DefContTheta} Definition We say that $\{\t_t\}_{t\in G}$ is
continuous if the corresponding map
  $$\t\:(a,t)\in \B^{-1} \arw (\t_t(a),t)\in \B$$ is continuous.

  Note that \scite{\Fell}{II.13.16} provides the following equivalent
characterization of continuity:

  \state Proposition Suppose that $\Gamma$ is a fixed pointwise-dense
space of sections for $\D^{-1}$. Then $\{\t_t\}_{t\in G}$ is
continuous if and only if for any $\gamma$ in $\Gamma$ one has that
the map
  $t\in G\mapsto \t_t(\gamma(t))\in A$
  is continuous.

  From \scite{\Fell}{II.13.17} it follows that, if $\t$ is continuous,
then its inverse is continuous as well.

  We must now discuss continuity for
  $\{w(r,s)\}_{(r,s)\in G\x G}$.
  The idea again will be to define continuity in terms of the
continuity of the corresponding bundle map. However, there is a slight
problem here because $w(r,s)$ is a map (actually a multiplier consists
of a pair of maps) defined in
  $D_r\cap D_{rs}$
  and one should worry
  in the first place whether or not these form a Banach bundle. The
question here is
  whether the pointwise intersection of two continuously varying
families of subspaces is again continuous in our sense. In general
this is not the case but, fortunately, this is true for ideals. In
fact, using \lcite{\LotsSections} and recalling that the intersection
of two ideals equals their product, one obtains enough sections of the
intersection bundle by multiplying together a pair of sections of each
bundle.
  So, let $\D$ be the Banach bundle over $G\x G$ having
  $D_r\cap D_{rs}$
  as the fiber over $(r,s)$ in the spirit of
\lcite{\BundleVsContinuity}. The family $\{w(r,s)\}_{(r,s)\in G\x G}$
then defines bundle maps
  $$L, R\:\D\arw \D$$
  given by the left and right action of the multipliers $w(r,s)$,
respectively.

  \nstate {\DefContW} Definition We say that $\{w(r,s)\}_{(r,s)\in G\x
G}$ is continuous if both $L$ and $R$ are continuous maps from $\D$ to
$\D$.  Equivalently (see \scite{\Fell}{II.13.16}), if for any $\gamma$
in a fixed pointwise-dense set of sections of $\D$ one has that both
  $$(r,s)\in G\mapsto \gamma(r,s)w(r,s)\in A$$ and
  $$(r,s)\in G\mapsto w(r,s)\gamma(r,s)\in A$$
  are continuous.

  The definition of continuity for \tpas is thus obtained by requiring
that all of its components be continuous in the appropriate senses:

  \state Definition If\/
 $\Theta = \(\{D_t\}_{t\in G}, \{\t_t\}_{t\in G}, \{w(r,s)\}_{(r,s)\in
G\x G}\)$
  is a \tpa of the \lc group $G$ on the \cstar-algebra $A$, we say
that $\Theta$ is continuous if
  \items
  \itm{a} $\{D_t\}_{t\in G}$ is continuous in the sense of
\lcite{\DefContSubspace}.
  \itm{b} $\{\t_t\}_{t\in G}$ is continuous in the sense of
\lcite{\DefContTheta}.
  \itm{c} $\{w(r,s)\}_{(r,s)\in G\x G}$ is continuous in the sense of
\lcite{\DefContW}.
  \endgroup

  In the special case in which the action is neither twisted nor
partial, that is, if $w(r,s)=1$ and $D_r=A$ for all $r$ and $s$ in
$G$, observe that our definition of continuity reduces to the usual
concept of strongly continuous group action.

  Let us now fix a continuous \tpa\null $\Theta$ of the \lc group $G$
on $A$. Our goal will be to show that the \spb of section \lcite{2},
equipped with the relative topology (from $G\x A$) is a continuous
\cstar-\ab over $G$. The above considerations, in particular
Proposition \lcite{\BundleVsContinuity}, already tells us that $\B$ is
a Banach bundle over $G$. For the convenience of the reader let us
recall the definition of a continuous \cstar-\ab.

  \state Definition A continuous \cstar-\ab over the \lc group $G$ is
a Banach bundle $\B=\{B_t\}_{t\in G}$ together with a continuous
multiplication operation
  $$\cdot\ \:\B\x \B \arw \B$$ and a continuous involution
  $$*\ \:\B\arw \B$$
  satisfying (i)--(xi) of \lcite{\CstarAlgBundle}.

  The following concludes the presentation of our main construction.

  \nstate{\Construction} Theorem Let
 $\Theta = \(\{D_t\}_{t\in G}, \{\t_t\}_{t\in G}, \{w(r,s)\}_{(r,s)\in
G\x G}\)$ be a continuous \tpa of the \lc group $G$ on the
\cstar-algebra $A$. Then the \spb of $A$ and $G$, with the relative
topology of $A\x G$, is a continuous \cstar-\ab over $G$.

  \beginproof After the work done in section \lcite{2}, it is enough
to verify the continuity of the multiplication and of the
involution. In order to do this, we use \scite{\Fell}{VIII.2.4 and
3.2}.
  Assume that $\alpha$ and $\beta$ are continuous sections of $\B$. So
$\alpha(t) = a_t\d_t$ and $\beta(t) = b_t\d_t$ where $a_t$ and $b_t$
are continuous $A$-valued functions on $G$ with $a_t,b_t\in D_t$ for
all $t$. We must therefore prove that the map
  $(r,s)\in G \mapsto \alpha(r)\beta(s)\in \B$ is continuous.
  By definition, we have
  $\alpha(r)\beta(s) = \t_r\(\t_r^{-1}(a_r)b_s\)w(r,s)\d_{rs}$. Now,
note that the map
  $$\gamma\:(r,s)\in G\x G \mapsto \t_r\(\t_r^{-1}(a_r)b_s\) \in A$$
  is continuous. Also, observe that, since $\gamma(r,s) \in D_r\cap
D_{rs}$, we see that $\gamma$ gives a continuous section of the bundle
$\D$, mentioned before \lcite{\DefContW}. The continuity of $w(r,s)$
can now be invoked to conclude that $\alpha(r)\beta(s)$ is continuous.

  With respect to the involution, using \scite{\Fell}{VIII.3.2}, we
must show that for each continuous section $\alpha(t) = a_t\d_t$, one
has that $\alpha(t)^*$ is continuous. Recall that $\alpha(t)^* =
\t_t^{-1}(a_t^*)w(t^{-1},t)^*\d_{t^{-1}}$. So proving continuity of
$\alpha(t)^*$ amounts to proving the continuity of the $A$-valued
function
  $$t\in G\mapsto \t_t^{-1}(a_t^*)w(t^{-1},t)^*\in A.$$

  Note that this map is given by the composition of the continuous map
  $$t\in G \mapsto \((t^{-1},t),\t_t^{-1}(a_t^*)\) \in \D$$
  followed by the inverse of the map $L$ mentioned in
\lcite{\DefContW}, and finally, the projection from $\D$ to $A$.
\endproof

  \beginsection Ternary Rings of Operators

  If $\B$ is a \cstar-\ab then, except for the unit fiber algebra
$B_e$, the fibers $B_t$ are not closed under multiplication and,
therefore, do not possess the structure of an algebra. Nevertheless
$B_t$ has a rich algebraic structure provided by the ternary operation
$xy^*z$, for $x,y,z\in B_t$, with respect to which it is closed.  This
makes $B_t$ a ternary \cstar-ring as defined by Zettl \cite{\Zettl}.

  \nstate{\DefESpace} Definition A ternary \cstar-ring is a complex
Banach space $E$, equipped with a ternary operation
  $$(a,b,c)\in E\x E\x E \mapsto \m abc \in E$$
  which is linear in the first and third variables, conjugate linear
in the middle variable and which satisfies the following for all
$a,b,c,d,e\in E$
  \items
  \itm{i} $\m {(\m abc)}de = \m a{(\m dcb )}e = \m ab{(\m cde)}$
  \itm{ii} $\|\m abc\| \leq \|a\|\,\|b\|\,\|c\|$
  \itm{iii} $\|\m aaa\| = \|a\|^3$
  \endgroup

  According to Theorem \lcite{3.1} in \cite{\Zettl}, there is a
fundamental dichotomy in the theory of ternary \cstar-rings in the
sense that any such object is the direct sum of a {\it ternary ring of
operators} (\tro) and what could be called an anti-\tro.  The
definitions are as follows:

  \nstate{\DefTRO} Definition A ternary ring of operators is a closed
subspace $E$ of operators on a Hilbert space $H$ (Zettl considers the
case of operators between two Hilbert spaces but this is not a crucial
matter) such that $EE^*E\c E$, equipped with the ternary operation
  $$\m abc = ab^*c \for a,b,c\in E.$$
  An anti-\tro is a \tro except that the ternary operation, with which
it comes equipped, is
  $$\m abc = -ab^*c \for a,b,c\in E.$$

  Clearly \tros as well as anti-\tros are examples of ternary
\cstar-rings. It is interesting to remark, however, that there is a
legitimate difference between these in the sense that a \tro is not
isomorphic to an anti-\tro or vice-versa. This is related to the
uniqueness in Theorem \lcite{3.1} of \cite{\Zettl}.

  Fortunately we will only have to deal with \tros here, mainly
because the fibers of a \cstar-\ab are \tros, a fact that follows from
\scite{\Fell}{VIII.16.5}.

  Loosely following \cite{\Zettl}, and occasionally offering minor
improvements, we propose to discuss below a few facts about \tros
which will be needed in the sequel. We will often use the notation
``$\m abc$'' in place of ``$ab^*c$'' because most of our results
concern the intrinsic structure of \tros, irrespective of the Hilbert
space representation which is, nevertheless, always in the background.

  Let $E$ be a \tro which we consider fixed for the time being.

  \state Definition A map $T\:E\arw E$ is said to be a left \rr
operator if there exists another map $T^*\:E\arw E$ satisfying
  $$\m a{T(b)}c = \m a b{T^*(c)}$$
  $$ \hbox{(resp.}\quad\m a{T(b)}c = \m {T^*(a)}bc)$$

  Note that, since the (ternary) multiplication on $E$ is
non-degenerate in view of \lcite{\DefESpace.iii}, $T^*$, if it exists,
must be unique. The definition also implies that a left \rr operator
must necessarily be a bounded linear map (for boundedness, use the
closed graph theorem).

  \state Proposition If\/ $T$ is a left \rr operator, then
  \items
  \itm{i} $T^*$ is also a left \rr operator and $T^{**}=T$.
  \itm{ii} For any $a,b,c \in E$ one has $\m {T(a)}bc = T(\m abc)$
(resp.  $\m ab{T(c)} = T(\m abc))$.
  \endgroup

  \beginproof Let us assume $T$ is a left operator and let
$a,b,c,x,y\in E$. Then
  $$\m x{(\m a{T^*(b)}c)}y = \m{(\m xc{T^*(b)})}ay$$
  $$ = \m{(\m x{T(c)}b)}ay = \m x{(\m ab{T(c)})}y$$
  This, together with the non-degeneracy of the product, implies that
$T^*$ is a left operator and that its adjoint is $T$.

  To prove \lcite{ii} we have
  $$\m x{(\m {T(a)}bc)}y = \m xc{(\m b{T(a)}y)}=\m xc{(\m
ba{T^*(y)})}$$
  $$=\m x {(\m abc)}{T^*(y)} = \m x{T(\m abc)}y.$$

  The proof for right operators is similar. \endproof

  \nstate{\LeftOperIsAlg} Proposition The set of all left \rr
operators on $E$ is a \cstar-algebra under the composition of
operators, the involution defined above and the operator norm.

  \beginproof Let $T$ be a left operator. For $x$ in $E$ we have
  $$\|T(x)\|^3 = \|\m {T(x)}{T(x)}{T(x)}\| = \|\m
{T(x)}{x}{T^*(T(x))}\| $$
  $$\leq \|T(x)\|\, \|x\|\, \|T^*(T(x))\| \leq \|T\|\, \|T^*T\|\,
\|x\|^3.$$
  This shows that $\|T\|^3 \leq \|T\|\, \|T^*T\|$ from which it
follows that $\|T\|^2 \leq \|T^*T\|$. This can now be used to show
both the norm preservation of the adjoint operation and the
\cstar-identity: $\|T\|^2 = \|T^*T\|$. The verification of the
remaining properties is left to the reader. \endproof

  The algebra of left operators on $E$ will be denoted $\L(E)$ and
will be called the left algebra. Similarly we have the right algebra
$\R(E)$. Please note that we are not using the same notation as in
\cite{\Zettl}.

  Given $x$ and $y$ in $E$, consider the maps $\lambda_{xy}\:E\arw E$
and $\rho_{xy}\:E\arw E$ defined by
  $$\lambda_{xy}(a) = \m xya \and \rho_{xy}(a) = \m axy$$
  for all $a$ in $E$.  Note that for $a,b,c \in E$ we have
  $$\m a{\lambda_{xy}(b)}c = \m a{(\m xyb)}c = \m ab{(\m yxc)} =
  \m ab{\lambda_{yx}(c)}$$
  so $\lambda_{xy}$ is a left operator and
$\lambda_{xy}^*=\lambda_{yx}$.  Similarly, $\rho_{xy}$ is a right
operator and $\rho_{xy}^*=\rho_{yx}$.

  If $T$ is a left operator, then one can easily show that
  $$T \lambda_{xy} = \lambda_{T(x),y}$$
  and that
  $$\lambda_{xy} T = \lambda_{x,T^*(y)}.$$
  So, one concludes that the closed linear span of the set of all
$\lambda_{xy}$, within $\L(E)$, is an ideal which we denote by $E\*
E^*$.  Similarly $E^*\* E$ is the ideal of $\R(E)$ given by the closed
linear span of the $\rho_{xy}$.

  Observe that, as a consequence of $EE^*E\c E$, one has that both
$EE^*$ and $E^*E$ are closed under composition of operators.

  Before we proceed, let us establish a slightly unusual notational
convention which will, nevertheless, serve our purposes rather well:

  \nstate{\Convention} Definition If $X$ and $Y$ are sets of elements
such that some kind of multiplication $xy$ is defined for $x$ in $X$
and $y$ in $Y$, taking values in some normed linear space, then $XY$
denotes the \underbar{closed linear span} of the set of products $xy$
with $x \in X$ and $y \in Y$. This applies, in particular, to subsets
of a \cstar-algebra and also when $X$ is a set of operators and $Y$ is
a set of vectors operated upon by the elements of $X$. The extreme
situation in which $X$ consists of a single element $1$, which acts as
a neutral element on $Y$, will be enforced as well. That is, $1Y$ is
the \underbar{closed linear span} of $Y$ rather than $Y$ itself.

  So $EE^*$ and $E^*E$, once interpreted according to the above
definition, are actually \cstar-algebras of operators on the Hilbert
space where $E$ acts. We would now like to prove that these are
isomorphic to $E\*E^*$ and $E^*\*E$, respectively. This should be
thought of as an indication that \tros are abstract objects which do
not depend so much on the representation considered. The precise
expression of this truth is Theorem \lcite{3.1} in \cite {\Zettl},
which we already mentioned.

  \nstate{\NormInE} Lemma Let $E$ be a \tro on the Hilbert space $H$.
If $a$ is in $EE^*$ and $b$ is in $E^*E$ then
  \items
  \itm{i} $\|a\| = \sup\{\|ax\|\: x\in E, \|x\| \leq 1\}$
  \itm{ii} $\|b\| = \sup\{\|xb\|\: x\in E, \|x\| \leq 1\}$
  \endgroup

  \beginproof We prove only \lcite{i}. Define a norm
  $\tri{\cdot}$ on $EE^*$ using the right hand side of \lcite{i}.
Since $E$ is invariant under left multiplication by $EE^*$, it follows
that $EE^*$ is a normed algebra under this new norm.

Note that for $c$ in $EE^*$ we have
  \def\sopa{\sup_{\|x\|= 1 \atop x\in E}}
  $$\tri{c}^2 = \sopa \|cx\|^2 = \sopa \|x^*c^*cx\| \leq \sopa
\|c^*cx\| = \tri{c^*c}$$
  Therefore $\tri{c}^2 \leq \tri{c^*c}$ which implies that $EE^*$ is a
pre-\cstar-algebra under this norm. But, since $\tri{a} \leq \|a\|$ we
must have $\tri{a} = \|a\|$.
  % Add reference here
  Part \lcite{ii} follows similarly.  \endproof

  \state Proposition Let $E$ be a \tro on $H$. Then there are
bijective \cstar-algebra isomorphisms
  $$\phi\:E\*E^* \arw EE^* \and \psi\: E^*\*E \arw E^*E$$
  such that
  $$\phi(\lambda_{xy}) = xy^* \and \psi(\rho_{xy}) = x^*y.$$

  \beginproof Let $\alpha\in E\*E^*$ be the finite sum, $\alpha =
\sum\lambda_{x_iy_i}$. Define $\phi(\alpha) = \sum x_i y_i^*$. To see
that this is well defined note that, with the help of Lemma
\lcite{\NormInE}, we have
  \def\sopa{\sup_{\|z\|= 1 \atop z\in E}}
  $$\|\sum x_iy_i^*\|
  = \sopa \| \sum x_iy_i^*z\| =$$
  $$ \sopa \| \sum \m {x_i}{y_i}z\|
  = \sopa \| \sum \lambda_{x_iy_i}(z)\| = \|\alpha\|.$$

  This shows that $\phi$ is well defined and isometric. The remaining
verifications are left to the reader. \endproof

  Note that the essential spaces of $EE^*$ and $E^*E$ are,
respectively, $EH$ and $E^*H$. In addition the members of $E$ should
be thought of as being operators from $E^*H$ to $EH$ since they all
vanish on the orthogonal complement of the former, and have their
image contained in the latter. The following fact, which we will use
frequently, has appeared in \scite{\Newpim}{Proposition 2.6}.

  \state Proposition If\/ $\{u_i\}_i$ is an approximate identity for
$EE^*$ (resp. $E^*E$), then for all $x$ in $E$ we have $\lim_i u_i x =
x$ (resp. $\lim_i x u_i = x)$.

  It is a consequence of this, that:

  \state Corollary If $E$ is a \tro then $EE^*E=E$.

  Let us now study the question of stability for \tros.

  \state Definition A \tro\null\ $E$ is said to be left \rr stable if
$E\*E^*$ (resp. $E^*\*E$) is a stable \cstar-algebra. In case $E$ is
both left and right stable, we simply say that $E$ is stable.

  A simple example of a \tro which is left-stable but not right-stable
is a Hilbert space equipped with the ternary multiplication
  $\m \xi \eta \zeta = \xi\<\eta,\zeta\>$ (where we think of
$\<\cdot\thinspace,\cdot\>$ as being conjugate-linear in the first
variable).

  In the following we let $\K$ denote the \cstar-algebra of all
compact operators on a separable infinite dimensional Hilbert space.

  \nstate{\KModuleIsStable} Proposition Suppose $E$ has the structure
of a left \rr module over $\K$ such that for all $a,b,c\in \K$ and $k$
in $\K$
  $$\m a{(kb)} c = \m ab{(k^*c)}$$
  $$\hbox{(resp.}\quad\m a{(bk)} c = \m {(ak^*)}bc).$$
  Suppose further that $E=\K E$ (resp. $E=E\K$). Then $E$ is left \rr
stable.

  \beginproof Suppose that $E$ is a left $\K$-module. Then the left
multiplication of elements of $K$ on $E$ gives a *-homomorphism $\K
\arw \L(E)$.  Recall that
  $$k\lambda_{xy} = \lambda_{kx,y}$$
  which implies that
  $\K (E\*E^*) \c E\*E^*$.
  Also, in view of the fact that $\K E=E$, and the equation above, we
actually conclude that
  $\K (E\*E^*) = E\*E^*$.

  It can be shown without much difficulty that if $A$ is a
\cstar-subalgebra of another \cstar-algebra, which also contains a
copy of the compact operators $\K$, such that $\K A = A$, then $A$
must be stable. Since this is precisely the case for $E\*E^*$, it
follows that this algebra is stable. The case of right $\K$-modules is
treated similarly. \endproof

  \beginsection Regular \tros\

  Given a \tro\null $E$, note that $E$ has a bi-module structure with
respect to the algebra $\L(E)$ acting on the left and $\R(E)$, on the
right (this means, in particular, that a left operator commutes with
any right operator \scite{\Zettl}{3.4}).
  In addition, $E$ is an imprimitivity bi-module for the ideals
$E\*E^*\c \L(E)$ and $E^*\*E\c \R(E)$
  \cite{\Induc}, \cite{\MoritaOne}, \cite{\MoritaTwo} with
inner-products defined by
  $$(x,y)\in E \arw (x|y) = \lambda_{xy}\in E\*E^*$$ and
  $$(x,y)\in E \arw \<x,y\> = \rho_{xy}\in E^*\*E.$$

The very delicate point of whether these inner products are positive
is not an issue here, precisely because we are dealing exclusively
with \tros, as opposed to general ternary \cstar-rings. In fact, under
the identification of $E\*E^*$ and $EE^*$, we have that $(x|x) = xx^*$
which is obviously a positive operator in $EE^*$. Likewise
$\<x,x\>=x^*x$ is positive.

  Closely associated with the notion of imprimitivity bi-modules,
there is the concept of linking algebra.  Recall that the linking
algebra \cite{\BGR}, \cite{\BMS} is the \cstar-algebra
  $$ \Link(E) = \mat{E\*E^*} {E} {E^*} {E^*\*E}$$
  equipped with the multiplication
  $$\mat{a_1} {x_1} {y_1^*} {b_1} \mat{a_2} {x_2} {y_2^*} {b_2}
  = \mat{a_1 a_2 + (x_1|y_2)} {a_1 x_2 + x_1 b_2}
    {y_1^* a_2 + b_1 y_2^*} {\<y_1, x_2\> + b_1 b_2}$$
  and involution
  $${\mat{a} {x} {y^*} {b}}^* = \mat{a^*} {y} {x^*} {b^*}$$
  for $a, a_1, a_2 \in E\*E^*$, $b, b_1, b_2 \in E^*\*E$ and $x, x_1,
x_2, y, y_1, y_2\in E$.

  Recall from \cite{\BGR} that two \cstar-algebras are said to be
strongly Morita equivalent to each other if there exists an
imprimitivity bi-module. Of course, whenever $E$ is a \tro, the
algebras $E\*E^*$ and $E^*\*E$ are Morita equivalent to each other.
If we add to this situation the hypothesis that $E$ is stable and that
both $E\*E^*$ and $E^*\*E$ possess strictly positive elements, then
the well known result of Brown, Green and Rieffel \scite{BGR}{Theorem
1.2} implies that $E\*E^*$ and $E^*\*E$ are isomorphic
\cstar-algebras.
  To pinpoint the precise consequence of this circle of ideas that we
will need is the main goal of the present section.

  In the following we let $E$ be a fixed \tro.

  \state Definition We say that $E$ is regular if there exists a
\piso\null $v$ in the multiplier algebra of\/ $\Link(E)$ such that
  $$vv^* = \mat{1}{0}{0}{0} \and v^*v=\mat{0}{0}{0}{1}.$$ (Compare
Lemma (3.3) of \cite{\BGR}).

  \nstate{\StableImpliesRegular} Proposition Any separable stable \tro
is regular.

  \beginproof This is an immediate consequence of combining (3.4) and
(3.3) in \cite{\BGR}, as long as we note that separability of $E$
implies separability of both $E\*E^*$ and $E^*\*E$ and hence the
existence of strictly positive elements. \endproof

  The following characterizes regular \tros at the Hilbert space
level.

  \nstate{\ConcreteU} Proposition Let $E$ be a \tro on $H$. Then $E$
is regular if and only if there exists a partially isometric operator
$u$ on $H$ such that
  $$uE^* = EE^* \and u^*E = E^*E.$$
  (Here as everywhere else in this work, we keep \lcite{\Convention}
in force. In particular $uE^*$, $u^*E$, $EE^*$ and $E^*E$ are all
meant to denote the closed linear span of the set of products).

  \beginproof Assume that $E$ is regular and hence assume the
existence of $v\in \M(\Link(E))$ as above. Denote by $e_1 = \mat 1000$
and $e_2=\mat 0001$, both of which are viewed as elements in
$\M(\Link(E))$. We have
  $$v^*\mat 0E00 = v^*e_1 \Link(E) e_2 = v^*vv^*\Link(E)e_2 =$$
  $$ e_2v^*\Link(E)e_2 \c e_2\Link(E)e_2
  = \mat 000{E^*\*E} = \mat 0E00^*\mat 0E00.$$
  Conversely
  $$ \mat 0E00^*\mat 0E00 = e_2\Link(E)e_2 = v^*vv^*v\Link(E)e_2 \c
v^*vv^* \Link(E)e_2$$
  $$=v^*e_1\Link(E)e_2 = v^*\mat 0E00.$$
  This proves that
  $$v^*\mat 0E00 = \mat 0E00^*\mat 0E00$$
  and in a similar way we could show that
  $$v\mat 0E00^* = \mat 0E00 \mat 0E00^*.$$

  Consider the representation of $\Link(E)$ on $H\+H$ given by
interpreting an element $\mat {a}{x}{y^*}{b}$ as an operator on $H\+H$
in the obvious way. Let $u$ be the image of $v$ under the canonic
extension of that representation to $\M(\Link(E))$. Since
  $vv^* = \mat{1}{0}{0}{0}$ and $v^*v=\mat{0}{0}{0}{1}$ it follows
that $u$ must actually have the form $\mat 0u00$, where $u\:H\arw H$
is a partially isometric operator.

  This said we see that the image of $v^*\mat 0E00$ in $\B(H\+H)$ will
thus be
  $$\mat00{u^*}0\mat 0E00 = \mat000{u^*E}$$
  which will coincide, by what we saw above, with $\mat
000{E^*E}$. This shows that $u^*E=E^*E$. Similarly $uE^*=EE^*$.

  To prove the converse statement, one defines $v$ to be the
multiplier on $\Link(E)$ whose left and right actions are given by
multiplying on the left and right by the operator $\mat
0u00$. Although
  $$\mat 0u00\mat 0u00^* = \mat 1000$$
  may not hold as operators on $H\+H$, that equality is true as long
as multipliers of $\Link(E)$ are concerned. Similarly $v^*v = \mat
0001$.\endproof

  \state Definition If\/ $E$ is a \tro and $u$ is a \piso such that
  $uE^* = EE^*$ and $u^*E = E^*E$, we say that $u$ is associated to
$E$ and we write $u\assoc E$.

  There is no obvious sense in which a \piso associated to $E$ is
unique.  In particular, the equations $uE^* = EE^*$ and $u^*E = E^*E$
do not even determine, in general, the initial and final space of
$u$. However we can at least affirm that the initial space of such a
$u$ contains $E^*H$. In fact
  $$E^*H = E^*EE^*H = u^*EE^*H \c u^*H.$$
  Likewise we have that $EH$ is contained in the image of $u$. Another
conclusion we can draw from the fact that $u\assoc E$, is that $u$
defines an isometry from $E^*H$ to $EH$. To see this note that
  $uE^*H = EE^*H \c EH$
  while
  $u^*EH = E^*EH \c E^*H$.

  \nstate{\DefStrict} Definition Let $u\assoc E$. Then $u$ is said to
be strictly associated to $E$ if the initial space of $u$ coincides
with $E^*H$. Equivalently, if the final space of $u$ is $EH$. In this
case we write $u\sassoc E$.

  Observe that, in case $u\assoc E$, but not strictly, then we can
make it strict by replacing $u$ by $up$ where $p$ is the orthogonal
projection onto $E^*H$. The property of being associated to $E$ will
not notice that change.

  A strict \piso also has a topological relationship to $E$:

  \nstate{\StrictlyIsInSClos} Proposition Let $E$ ba a \tro on $H$ and
assume that the partial isometry $u$ is associated to $E$. Then a
necessary and sufficient condition for $u$ to be strictly associated
to $E$ is that $u$ be in the strong operator closure of $E$ within
$\B(H)$.

  \beginproof
  Let $\{e_i\}_i$ be an approximate identity for $EE^*$. Then it is
well known that $e_i$ converges strongly to the the orthogonal
projection onto the essential space of $EE^*$, which we have seen to
coincide with $EH$. Let $u_i = e_i u$. Then $u_i \in EE^*u = EE^*E =
E$. Observing that the range of $u$ is $EH$, we conclude that $u_i$
converges strongly to $u$.  Conversely, if $u$ is in the strong
operator closure of $E$, then it must vanish on the orthogonal
complement of $E^*H$, as is the case for any member of $E$.  This
concludes the proof.  \endproof

  Returning to our earlier discussion on the question of uniqueness
for a \piso associated to $E$ observe that even strict ones are not
unique. In fact if $u\sassoc E$ and if $w$ is a unitary multiplier of
$EE^*$, then $wu\sassoc E$.

  With respect to the nature of the product $wu$ above, we need to
clarify the following point. If $A$ is a \cstar-algebra which is
represented on a Hilbert space $H$ under a faithful {\it
non-degenerate} representation, then it is well known
\scite{\Pedersen}{3.12.3} that its multiplier algebra $\M(A)$ can also
be represented within $\B(H)$.  However, if that representation is not
non-degenerated, i.e, if the essential space $AH$ is a proper subspace
of $H$, then this is still true in the sense that $\M(A)$ is
isomorphic to the algebra consisting of those operators $m$ in $\B(H)$
such that both $mA$ and $Am$ are contained in $A$ and such that both
$m$ and $m^*$ vanish on the orthogonal complement of $AH$. Therefore,
when we spoke of $w$ above, we meant an operator on $H$ and hence $wu$
should be simply interpreted as the composition of operators.

  \nstate{\UniqueUpToMultiplier} Proposition Let $u_1$ and $u_2$ be
strict \pisos associated to $E$. Then $u_2u_1^*$ is a unitary
multiplier in $\M(EE^*)$ and $u_2^*u_1$ is a unitary multiplier in
$\M(E^*E)$.

  \beginproof Proving the first statement amounts to verifying that
$u_2u_1^*$ is a unitary operator on $EH$ and that $EE^*$ is invariant
under both left and right multiplication by $u_2u_1^*$, all of which
follow by routine arguments. \endproof

  Still under the notation above, note that if $w=u_2u_1^*$ , then
$u_2=wu_1$ and so we see that strict \pisos are unique, after all, up
to multiplication by a unitary element in $\M(EE^*)$.

  \beginsection Ideals of \tros\

  The classification of \cstar-\abs we are about to discuss requires a
careful understanding of the relationship between \tros and its
subspaces, specially when these are ideals in the sense below.

  \nstate{\DefIdeal} Definition Let $J$ be a closed subspace of the
\tro\null $E$. Then $J$ is said to be an ideal if
  $$\m JJE \c J \and \m EJJ \c J.$$

We remark that there is a total number of eight possible ways of
combining $E$ and $J$ in the ternary product, so there are many
alternatives to the definition of the concept of ideals in \tros. Even
though we don't claim to have experimented with too many of those, we
hope to convince the reader that our choice is meaningful.

  Note that $\m JJE$ above, means $JJ^*E$ and a similar remark applies
to $\m EJJ$. From now on we will use the latter notation.

  \nstate{\Aux} Lemma If $J$ is an ideal in $E$, then
  \items
  \itm{i} $JJ^*E = J$
  \itm{ii} $EJ^*J = J$
  \itm{iii} $JJ^*EE^* = JJ^*$.
  \itm{iv} $J^*JE^*E = J^*J$.
  \endgroup

  \beginproof Initially observe that $J$ is a \tro in its own right.
With respect to \lcite{i} we have
  $$J = JJ^*J \c JJ^*E \c J$$
  so $J =JJ^*E$. The proof of \lcite{ii} is similar. As for
\lcite{iii} we have, using \lcite{i} and \lcite{ii}
  $$JJ^*EE^* = JE^* = JJ^*JE^* = J(EJ^*J)^* = JJ^*.$$
  We leave \lcite{iv} to the reader. \endproof

  \nstate{\RegularIdeal} Proposition An ideal $J$ of a regular
\tro\null $E$ is necessarily regular.  In addition if $u$ is a \piso
with $u\assoc E$, then $u\assoc J$.

  \beginproof It obviously suffices to prove the second assertion. For
that purpose we use \lcite{\Aux} in the following calculations
  $$uJ^* = uE^*JJ^* = EE^*JJ^* = JJ^*$$
  and
  $$u^*J = u^*EJ^*J = E^*EJ^*J = J^*J. \endproof$$

  Another fact we need for future use is proven below.

  \nstate{\ProdESpaces} Proposition Let $E$ and $F$ be \tros on $H$,
such that $FF^*E^*E = E^*EFF^*$, then
  \items
  \itm{i} $EF$ is a \tro.
  \itm{ii} If $u\sassoc E$ and $v\sassoc F$ then $uv\sassoc EF$.
  \endgroup

  \beginproof That $EF$ is a \tro follows from the following
calculation
  $$EFF^*E^*EF = EE^*EFF^*F =EF.$$

  To prove \lcite{ii} we first claim that the final projection $vv^*$
of $v$ commutes with the initial projection $u^*u$ of $u$. To see this
it is enough to show that the range of $vv^*$, which coincides with
$FH$, is invariant under $u^*u$. With that goal in mind note that $uFH
\c EFH$ because $u$ is in the strong closure of $E$. Hence
  $$u^*u(FH) \c u^*EFH = E^*EFF^*FH = FF^*E^*EFH \c FH.$$
  This shows our claim that $u^*u$ and $vv^*$ commute and hence that
$uv$ is a \piso.  To show that $uv \assoc EF$ we compute
  $$uvF^*E^* = uFF^*E^*EE^* = u E^*EFF^*E^* = EFF^*E^*.$$
  That $(uv)^*EF = (EF)^*EF$ follows similarly. It now remains to show
that the final space of $uv$ is $EFH$. We have
  $$uvH = uFH \c EFH$$
  again because $u$ is in the strong closure of $F$. The opposite
inclusion follows from the argument preceding \lcite{\DefStrict}.
\endproof

  Our motivation for conducting a study of \tros is, of course, our
interest in \cstar-\abs. If $\B$ is a \cstar-\ab and if $r$ and $s$
are elements in the base group, then it can be proven that $B_rB_s$ is
an ideal in $B_{rs}$. If everything is regular and represented in a
Hilbert space, then we will have isometries $u_r$, $u_s$ and
$u_{rs}$. To understand the relationship between $u_ru_s$ and $u_{rs}$
is our last objective before we plunge into the main section of this
work.

  \nstate{\ToGetMultiplier} Proposition Suppose $E$, $F$ and $M$ are
regular \tros such that
  \items
  \itm{i} $FF^*E^*E = E^*EF^*F$
  \itm{ii} $EF$ is an ideal in $M.$
  \endgroup

\noindent {\sl Let $u$, $v$ and $z$ be \pisos strictly associated to
$E$, $F$ and $M$, respectively. Then $uvz^*$ is a unitary multiplier
of $EFF^*E^*$.}

  \beginproof We know from \lcite{\ProdESpaces} that $uv\sassoc EF$.
On the other hand \lcite{\RegularIdeal} tells us that also $z\assoc
EF$ although possibly not strictly. But if $p$ is the orthogonal
projection onto $(EF)^*H$, then $zp\sassoc EF$. It then follows from
\lcite{\UniqueUpToMultiplier} that the operator $w$ defined by
$w=uv(zp)^*$ is a unitary multiplier in $\M(EFF^*E^*$). Finally
observe that since $uv$ is strict, then we must have $uvp=uv$ and
hence $w=uvz^*$. \endproof

  \beginsection The Classification of Stable \cstar-Algebraic Bundles

  Let $\B$ be a \cstar-\ab over the \lc group $G$, considered fixed
throughout this section. Our goal is to show, upon assuming a certain
regularity property of $\B$, that it is isomorphic to the \spb
constructed from a suitable \tpa of the base group $G$ on the unit
fiber algebra $B_e$.  Construction of that action will be done in
several steps.

  Let us initially deal with the problem of defining the family
$\{D_t\}_{t\in G}$ of ideals of $B_e$. We simply let, for each $t$ in
$G$, $D_t = B_tB_t^*$. Clearly $D_t$ is an ideal in $B_e$. To see that
the $D_t$ form a continuous family, recall from \lcite{\LotsSections}
that this follows once we provide a pointwise-dense set of sections.
Now observe that if $\gamma$ and $\delta$ are in $\cob$, the space of
continuous sections of $\B$ vanishing at infinity
\scite{\Fell}{II.14.7}, then $\gamma(t)\delta(t)^*$ is a continuous
$B_e$ valued function which satisfies $\gamma(t)\delta(t)^*\in D_t$
for all $t$. In other words it is a continuous section for the family
$\{D_t\}_{t\in G}$. The linear span of the set of such sections is
clearly pointwise-dense and hence this proves continuity.

  So far we have thus been able to construct a Banach bundle over $G$
  {}from the $D_t$'s, according to \lcite{\BundleVsContinuity}. Let us
denote the total space of this bundle by $\D$.

  Recall from \scite{\Fell}{II.14.1} that $\cob$ is a Banach space
under the supremum norm. It matters to us that it is also a \tro under
the ternary operation $(\m \g \d \e )(t) =
\g(t)\d(t)^*\e(t)$. Precisely, we mean that $\cob$ is isomorphic to a
\tro in some Hilbert space, as far as its Banach space structure and
the operation mentioned above are concerned. Let us now show how to
represent $\cob$ as a \tro.

  Let $\rho$ be a representation of $\B$ in the sense of
\scite{\Fell}{VIII.8.2, VIII.9.1} such that the restriction of $\rho$
to each $B_t$ is isometric. The existence of such a representation
follows from \scite{\Fell}{VIII.16.5}. We may therefore assume that
each $B_t$ is a closed subspace of $\B(H)$ for some Hilbert space $H$
(which does not depend on $t$) and moreover we have $B_rB_s \c B_{rs}$
and $B_t^* = B_{t^{-1}}$. In particular this implies that each $B_t$
is a \tro. Denoting by $\ltwo(G)$ the Hilbert space of all square
summable sequences of complex numbers, indexed by $G$ (we are
temporarily ignoring the topology of $G$ here), consider the map
  $$\pi\: \cob \arw \B\(H\*\ltwo(G)\)$$
  given by $\pi(\g)(\xi\*e_t) = (\g(t)\xi)\*e_t$
  for all $\g$ in $\cob$, $\xi$ in $H$ and $t$ in $G$, where we are
denoting by $e_t$ the canonic basis of $\ltwo(G)$. In other words
$\pi(\g)$ is the diagonal operator diag$(\g(t)_{t\in G})$, with
respect to the decomposition of $H\*\ltwo(G)$ provided by the canonic
basis of $\ltwo(G)$. Clearly $\pi$ is an isometric representation of
$\cob$ which satisfies
  $$\pi(\m \g \d \e) = \pi(\g) \pi(\d)^* \pi(\e).$$
  In other words we may identify $\cob$ with its image in
$\B\(H\*\ltwo(G)\)$ through $\pi$. This identification will be tacitly
made henceforth, without explicit mention to $\pi$.

In order to be able to apply the machinery of regular \tros developed
above, we would like to have $\cob$ regular. However this is not
within our reach unless we make extra requirements.

  \state Proposition If\/ $B_e$ is a stable \cstar-algebra then $\cob$
is stable as a \tro. If, in addition, $\B$ is second countable, then
$\cob$ is regular.

  \beginproof There is an obvious way in which $\cob$ can be
considered as a $B_e$ bi-module. With a little more effort we can give
$\cob$ the structure of a bi-module over the multiplier algebra
$\M(B_e)$. In order to do this we use \scite{\Fell}{VIII.3.8 and
VIII.16.3} to identify the multiplier algebra of $B_e$ with the set of
multipliers of $\B$ of order $e$ (see \scite{\Fell}{VIII.2.14}). Thus,
if $\g$ is in $\cob$ and $m$ is in $\M(B_e)$, we define $m\g$ in
$\cob$ by
  $(m\g)(t) = m (\g(t))$ and likewise
  $(\g m)(t) = (\g(t))m$.

  That $m\g$ and $\g m$ are continuous follows from the fact that the
left and right action of $m$ is a continuous bundle map from $\B$ to
itself. The latter, in turn, follows from \scite{\Fell}{II.13.16} with
the set $\Gamma =\{b_1\g b_2\in \cob\: b_1,b_2\in B_e, \g\in \cob\}$.
The reason why $\Gamma$ is pointwise dense, finally follows from the
existence of approximate identities \scite{\Fell}{VIII.16.3}.

  Observe that these module structures are compatible with the
*-operation in the sense that
 $$\m \g {m\d} \e = \m \g\d{m^*\e} \and
 \m \g {\d m} \e = \m {\g m^*}\d\e.$$
  In other words we have *-homomorphisms
  $$\Lambda\:\M(B_e) \arw \L(\cob))$$
  and
  $$\Rho\:\M(B_e) \arw \R(\cob).$$

  Given that $B_e$ is stable, say $B_e = A\*\K$ for some
\cstar-algebra $A$, consider the left and right actions of $K$ on
$B_e$ given by
  $k_1(a\*k_2):= a\*(k_1k_2)$
  and
  $(a\*k_2)k_1:= a\*(k_2k_1)$ for $k_1,k_2\in \K$ and $a\in A$. These
satisfy $\K B_e = B_e\K = \B_e$ and define a *-homomorphism $\phi\:\K
\arw \M(B_e)$. If we now follow this map by either $\Lambda$ or $\Rho$
above, we will give $\cob$ the structure of a $\K$ bi-module.
 We must now prove that $\K\cob = \cob \K = \cob$ to be in condition
to apply \lcite{\KModuleIsStable}. This follows from
  $$\K\cob = \K B_e\cob = B_e\cob = \cob$$
  and the corresponding right hand sided version.

  As for the second assertion in the statement, assume that $\B$ is
second countable. Then \scite{\Fell}{II.14.10} tells us that the space
of compactly supported continuous sections is separable in the
inductive limit topology. From this we can then deduce that $\cob$ is
separable in the sup norm. The conclusion is thus reached, upon
invoking \lcite{\StableImpliesRegular}. \endproof

  Although we will not explicitly need it, stability of $B_e$ implies
stability of each fiber $B_t$ as well.

  Let us assume, from now on, that $\B$ is a \cstar-\ab for which
$\cob$ is a regular \tro. Of course this will include all second
countable \cstar-\abs for which the unit fiber algebra is stable.

  As mentioned earlier in this section, for each pair of sections $\g$
and $\d$ in $\cob$, we have that $\g\d^*$ is a section of $\cod$. Note
that $\cod$ can also be represented in $H\*\ltwo(G)$ via diagonal
operators. Under these representations we therefore have
  $\cob\cob^* \c \cod$.  Using \scite{\Fell}{II.14.7} we actually
obtain
  $\cob\cob^* = \cod$.

  Recall that $\D^{-1}$ denotes the Banach bundle over $G$ which is
obtained by placing the ideal $D_{t^{-1}}$ as the fiber over $t$,
according to \lcite{\BundleVsContinuity}. Representing $\codinv$ also
via diagonal operators (each $D_{t^{-1}}$ acting on
  $H\*e_t$) we will have, by a similar reasoning, that $\cob^*\cob =
\codinv$.

  Since we are assuming that $\cob$ is regular, we may invoke
\lcite{\ConcreteU} to conclude that there must exist a \piso\null $U$
in $\B\(H\*\ltwo(G)\)$ which is strictly associated with $\cob$. By
\lcite{\StrictlyIsInSClos} we conclude that $U$ is in the strong
closure of $\cob$ within $\B\(H\*\ltwo(G)\)$ which, in turn, implies
that $U$ must be a diagonal operator. That is,
$U=\hbox{diag}\((u_t)_{t\in G}\)$, where, for for each $t$ in $G$,
$u_t$ is a \piso in $\B(H)$. Upon replacing $u_t$ by $u_e^*u_t$ we may
assume that $u_e=1$.

  Expressing in formulas the fact that $U\sassoc \cob$ we have
  $$U\cob^* = \cob\cob^* = \cod$$
  and
  $$U^*\cob = \cob^*\cob = \codinv.$$
  Under point evaluation at each group element $t$ (meaning to focus
on a specific diagonal entry) we see that each $u_t\assoc B_t$.
Moreover, since $U$ is strict with respect to $\cob$ we can easily
prove that $u_t$ is strict with respect to $B_t$.

  Given that $D_t = B_tB_t^* = B_t u_t^*$, we may define, for each $t$
in $G$, a map
  $$b\in B_t \arw bu_t^*\in D_t$$
  which is clearly an isometry onto $D_t$ and hence provides a bundle
map
  $$\rhodag\:\B \arw \D.$$
  That $\rhodag$ is continuous follow from \scite{\Fell}{II.13.16} and
the remark that for any continuous section $\g$ of $\B$ one has that
  $$\rhodag(\g(t)) = \g(t)u_t^* = (\g U^*)(t)$$
  which is a member of $\cob\cob^*=\cod$ and hence is continuous.  It
now follows from \scite{\Fell}{II.13.17} that $\rhodag$ is an
isometric isomorphism of Banach bundles. Its inverse is clearly the
map
  $$\rho\:\D \arw \B$$
  given by gluing together the maps
  $$a\in D_t \arw a u_t\in B_t.$$

  In an entirely similar way we have the isometric Banach bundle
isomorphism
  $$\lamdag\:\B \arw \D^{-1}$$
  which, together with its inverse
  $$\lambda\:\D^{-1} \arw \B$$
  are given by
  $$\lamdag\:b\in B_t \arw u_t^*b\in D_{t^{-1}}
  \and
  \lambda\:b\in D_{t^{-1}} \arw u_tb\in B_t.$$

  Now define, for each $x$ in $D_{t^{-1}}$,
  $$\t_t(x)=u_txu_t^*.$$
  Note that
  $$\t_t(D_{t^{-1}}) = \t_t(B_t^*B_t) = u_tB_t^*B_tu_t^* = B_tB_t^* =
D_t.$$
  Since $u_t^*u_t$ is the identity on $B_t^*H$, we see that $\t_t$ is
a \cstar-algebra isomorphism. The continuity of $\t$, in the sense of
\lcite{\DefContTheta}, is now obvious since the corresponding bundle
map $\theta\:\D^{-1} \arw \D$ is just the composition
  $$\D^{-1} \stackrel{\lambda}{\arw} \B \stackrel{\rho^\dagger}{\arw}
\D.$$

  The only missing ingredient of our \tpa is now the multipliers
$w(r,s)$ mentioned in \lcite{\DefTPA}. Referring to the notation of
\lcite{\ToGetMultiplier}, let, for each $r$ and $s$ in $G$, $E=B_r$,
$F=B_s$ and $M=B_{rs}$. Once one verifies that the hypothesis of
\lcite{\ToGetMultiplier} are verified, we will conclude that the
element $w(r,s)$ defined by $w(r,s)=u_ru_su_{rs}^*$ is a unitary
multiplier of $B_rB_sB_s^*B_r^*$. As for
  \lcite{\ToGetMultiplier.i}, this holds because both $B_sB_s^*$ and
$B_r^*B_r$ are ideals in $B_e$ and hence their product in either order
coincides with their intersection.

  \state Lemma For all $r$ and $s$ in $G$ one has $B_rB_sB_s^*B_r^*=
D_r\cap D_{rs}$.

  \beginproof We have
  $$B_rB_sB_s^*B_r^*
  = B_rB_r^*B_rB_sB_s^*B_r^*
  \c B_rB_r^*B_{rs}B_{rs}^*
  = D_r\cap D_{rs}.$$
  Conversely
  $$D_r\cap D_{rs} =
  D_r D_{rs} =
  D_r D_{rs} D_r =
  B_rB_r^*B_{rs}B_{rs}^*B_rB_r^* \c
  B_rB_sB_s^*B_r^*.\endproof$$

  With this at hand we see that $w(r,s)$ is a multiplier of $D_r\cap
D_{rs}$ as it is called for by \lcite{\DefTPA}. To see that $w(r,s)$
is continuous as a section of the bundle formed by the $D_r\cap
D_{rs}$ we can employ \scite{\Fell}{II.14.16} using the
pointwise-dense space of sections that is spanned by the sections of
the form
  $(r,s)\mapsto \g(r)\d(s)\e(s)^*\zeta(r)^*$ where $\g,\d,\e$ and
$\zeta$ are in $\cob$.

  The following is our main result:

  \nstate{\Main} Theorem Let $G$ be a \lc group and let $\B$ be a
\cstar-\ab over $G$ such that $\cob$ is regular (e.g. if $\B$ is
second countable and $B_e$ is stable). Then, there exists a continuous
twisted partial action
  $$\Theta = \(\{D_t\}_{t\in G}, \{\t_t\}_{t\in G}, \{w(r,s)\}_{(r,s)
\in G\x G}\)$$
  of $G$ on the unit fiber algebra $B_e$, such that $\B$ is
isometricaly isomorphic to the \spb of $A$ and $G$ constructed from
$\Theta$.

  \beginproof Define $D_t$, $\theta_t$ and $w(r,s)$ as above, and let
us prove that
  $$\Theta = \(\{D_t\}_{t\in G}, \{\t_t\}_{t\in G}, \{w(r,s)\}_{(r,s)
\in G\x G}\)$$
  is a continuous \tpa. Starting with \lcite{\DefTPA.b} we have
  $$\t_r(D_{r^{-1}}\cap D_s) = \t_r(B_r^*B_rB_sB_s^*) =
  u_rB_r^*B_rB_sB_s^*B_r^*B_ru_r^*=$$
  $$B_rB_sB_s^*B_r^* = D_r\cap D_{rs}.$$

  To prove \lcite{\DefTPA.c} let $a$ be in $D_{s^{-1}}\cap
D_{s^{-1}r^{-1}}$. Then
  $$w(r,s)\t_{rs}(a)w(r,s)^* =
  u_r u_s u_{rs}^* u_{rs} a u_{rs}^* u_{rs} u_s^* u_r^*.$$
  Now, since $u_{rs}^* u_{rs}$ is the projection on $E_{rs}^*H$, it
follows that $u_{rs}^* u_{rs}$ behaves like the identity element when
it is multiplied by elements in $E_{rs}^*E_{rs} = D_{s^{-1}r^{-1}}$.
This proves \lcite{\DefTPA.c}.

  With respect to \lcite{\DefTPA.d}, recall that $u_e=1$ so that
  $$w(t,e)=u_tu_eu_t^* = u_tu_t^*$$
  which is precisely the unit in the multiplier algebra $\M(D_t)$, at
least according to our convention discussed before
\lcite{\UniqueUpToMultiplier}.
  Similarly $w(e,t)=1$.

  The last axiom in \lcite{\DefTPA} translates to the following, for
  $a \in D_{r^{-1}}\cap D_s\cap D_{st}$
  $$u_rau_s u_t u_{st}^* u_r^*u_r u_{st} u_{rst}^* =
  u_r a u_r^* u_r u_s u_{rs}^* u_{rs} u_t u_{rst}^*.$$
  To see that this holds, all we need is to show that the initial
projections of the various partial isometries appearing in this
expression may be canceled out. This can be done by observing that, in
all cases, that projection appears besides an operator which `lives'
in its range.

  Now, given that the continuity of the various ingredients of our
\tpa have already been verified, we conclude that $\Theta$ is indeed a
continuous \tpa.

  Let, therefore, $\D$ be the \cstar\-\ab obtained from $\Theta$ as
described in \lcite{\Construction}. To conclude we must prove that
$\B$ and $\D$ are isomorphic \cstar-\ab. Recall that we have already
found a isometric Banach bundle isomorphism
  $$\rho\:\D \arw \B$$
  given by
  $$\rho(a_t\d_t) = a_tu_t, \for a_t \in D_t,$$
  which we now claim to be a \cstar-\ab isomorphism as well. To prove
this claim all we need to check is that $\rho$ is multiplicative and
*-preserving. For $a_r$ in $D_r$ and $b_s$ in $D_s$ let us prove that
  $$\rho(a\d_r * b_s\d_s) = \rho(a_r\d_r) \rho(b_s\d_s).$$
  The left hand side equals
  $$\rho(\t_r\(\t_r^{-1}(a_r)b_s\)w(r,s)\d_{rs}) =
  \t_r\(\t_r^{-1}(a_r)b_s\)u_ru_su_{rs}^*u_{rs} =$$
  $$ u_r(u_r^* a_r u_r b_s)u_r^* u_ru_su_{rs}^*u_{rs} =
  a_r u_r b_s u_s = \rho(a_r\d_r) \rho(a_s\d_s).$$

  Finally
  $$\rho\((a_r\d_r)^*\)
  = \rho\(\t_r^{-1}(a_r^*)w(r^{-1},r)^*\d_{r^{-1}}\)
  = u_r^*a_r^* u_r u_r^* u_{r^{-1}}^* u_{r^{-1}} = u_r^* a_r^* =
(a_ru_r)^*. \endproof$$

  \bigbreak
  \centerline{\tensc References}
  \nobreak\medskip
  \frenchspacing

  \stdbib{\BGR}
  {L. G. Brown, P. Green and M. A. Rieffel}
  {Stable isomorphism and strong Morita equivalence of
\cstar-algebras}
  {Pacific J. Math.} {71} {1977} {349} {363}

  \bib{\BMS}
  {L. G. Brown, J. A. Mingo and N. T. Shen}
  {Quasi-multipliers and embeddings of Hilbert \cstar-bimodules}
  {preprint, Queen's University, 1992}

  \bib{\Newpim}
  {R. Exel}
  {Circle actions on \cstar-algebras, partial automorphisms and a
  generalized Pimsner--Voiculescu exact sequence}
  {{\sl J. Funct. Analysis}, to appear}

  \bib{\Fell}
  {J. M. G. Fell and R. S. Doran}
  {Representations of *-algebras, locally compact groups, and Banach
*-algebraic bundles}
  {Academic Press, Pure and Applied Mathematics, vols. 125 and 126,
1988}

  \bib{\McCl}
  {K. McClanahan}
  {$K$-theory for partial crossed products by discrete groups}
  {preprint, University of Mississipi}

  \stdbib{\PR}
  {J. A. Packer and I. Raeburn}
  {Twisted products of \cstar-algebras}
  {Math. Proc. Camb. Phil. Soc} {106} {1989} {293} {311}

  \bib{\Pedersen}
  {G. K. Pedersen}
  {\cstar-Algebras and their Automorphism Groups}
  {Academic Press, 1979}

  \bib{\Quigg}
  {J. C. Quigg}
  {Discrete \cstar-coactions and \cstar-\ab}
  {preprint, Arizona State University}

  \stdbib{\Induc}
  {M. A. Rieffel}
  {Induced representations of \cstar-algebras}
  {Adv. Math.} {13} {1974} {176} {257}

  \stdbib{\MoritaOne}
  {M. A. Rieffel}
  {Morita equivalence for \cstar-algebras and $W^*$-algebras}
  {J. Pure Appl. Algebra} {5} {1974} {51} {96}

  \stdbib{\MoritaTwo}
  {M. A. Rieffel}
  {Strong Morita equivalence of certain transformation group
\cstar-alge\-bras}
  {Math. Ann.} {222} {1976} {7} {22}

  \stdbib{\ZM}
  {G. Zeller-Meier}
  {Produits crois\'es d'une \cstar-alg\`ebre par un group
d'automorphismes}
  {J. Math Pures Appl.} {47} {1968} {101} {239}

  \stdbib{\Zettl}
  {H. Zettl}
  {A characterization of ternary rings of operators}
  {Adv. Math.} {48} {1983} {117} {143}

  \vskip 2cm
  \rightline{May 2, 1994}

  \bye